\documentclass[journal]{IEEEtran}
\usepackage{cite}
\usepackage{amsmath,amssymb,amsfonts}
\usepackage{algorithm}
\usepackage{algpseudocode}
\usepackage{graphicx} 
\usepackage{xcolor}
\usepackage{siunitx}
\def\BibTeX{{\rm B\kern-.05em{\sc i\kern-.025em b}\kern-.08em
    T\kern-.1667em\lower.7ex\hbox{E}\kern-.125emX}}
\usepackage{acro}
\usepackage{stfloats}
\usepackage{caption}
\usepackage{subcaption}
\usepackage{url}
\usepackage{hyperref}
\usepackage{tikz,pgfplots}
\usepackage{comment}
\usepackage{amsthm}
\theoremstyle{remark}
\newtheorem*{remark}{Remark}

\usepackage{textcomp}

\usetikzlibrary{positioning}
\usetikzlibrary{calc}
\usetikzlibrary{arrows}
\usetikzlibrary{shapes.geometric}

\newcommand{\dc}{\mathrm{dc}}
\newcommand{\abc}{\mathrm{abc}}
\newcommand{\A}{\mathrm{a}}
\newcommand{\B}{\mathrm{b}}
\newcommand{\C}{\mathrm{c}}

\newcommand{\s}{\mathrm{s}}
\newcommand{\m}{\mathrm{m}}
\newcommand{\e}{\mathrm{e}}
\newcommand{\R}{\mathbb{R}}
\newcommand{\SM}{\mathbb{S}}
\newcommand{\factor}{\mathrm{fct}}
\newcommand{\p}{\mathrm{p}}
\newcommand{\diag}{\mathrm{diag}}
\newcommand{\U}{\mathrm{u}}
\newcommand{\opt}{\mathrm{opt}}
\newcommand{\T}{\mathrm{T}}
\newcommand{\ini}{\mathrm{ini}}
\newcommand{\inc}{\mathrm{inc}}
\newcommand{\prev}{\mathrm{prev}}
\newcommand{\trace}{\mathrm{tr}}

\newcommand{\sw}{\mathrm{sw}}
\newcommand{\rel}{\mathrm{sdp}}

\newcommand{\nodelim}{500}

\renewcommand{\r}{\mathrm{r}}
\renewcommand{\b}[1]{\boldsymbol{#1}}

\newlength{\fheight}
\newlength{\fwidth}

\DeclareSIUnit{\pu}{pu}
\DeclareSIUnit{\voltampere}{VA}

\DeclareAcronym{3L-NPC}{short = 3L-NPC, long = three-level neutral-point-clamped}
\DeclareAcronym{c1}{short = DTFCS-MPC, long = direct-torque FCS-MPC}
\DeclareAcronym{MPC}{short = MPC, long = model predictive control}
\DeclareAcronym{FCS-MPC}{short = FCS-MPC, long = finite-control-set MPC}
\DeclareAcronym{SDP}{short = SDP, long = semidefinite programming}
\DeclareAcronym{PSD}{short = PSD, long = positive semidefinite}
\DeclareAcronym{ED}{short = EVD, long = eigenvalue decomposition}
\DeclareAcronym{IM}{short = IM, long = induction motor}    
\begin{document}

\title{On the SDP Relaxation of Direct Torque Finite Control Set Model Predictive Control}

\author{Luca M. Hartmann, Orcun Karaca, Tinus Dorfling, Tobias Geyer, Adam Kurpisz
\thanks{Corresponding author: O. Karaca.\\ \indent L. M. Hartmann is with University of California San Diego, USA. email: {\tt lhartmann@ucsd.edu.} O. Karaca and T. Dorfling are with ABB Corporate Research, Switzerland. emails: {\tt \{orcun.karaca,} {\tt martinus-david.dorfling\}@ch.abb.com}. T. Geyer is with ABB System Drives, Switzerland. email: {\tt t.geyer@ieee.org}. A. Kurpisz is with ETH Zürich and BFH, Switzerland. email: {\tt adam.kurpisz@ifor.math.ethz.ch}.}}

\maketitle

\begin{abstract}
This paper formulates a semidefinite programming relaxation for a long horizon direct-torque finite-control-set model predictive control problem. 
In parallel with this relaxation, a conventional branch-and-bound algorithm tailored for the original problem, but with an iteration limit to restrict its computational burden, is also solved. 
An input sequence candidate is extracted from the solution of the semidefinite program in the lifted space. This sequence is then compared with the so-called early-stopping branch-and-bound solution, and the best of the two is applied in a receding horizon fashion. In simulated case studies, the proposed approach exhibits significant improvements in torque transients, as the branch-and-bound alone struggles to find a meaningful solution due to the imposed limit.
\end{abstract}


\section{Introduction}

Direct \ac{MPC} methods, which combine control and modulation into one stage by directly manipulating the switching positions of a converter, have garnered significant attention within the power electronics community~\cite{karamanakos2020model,rodriguez2021latest,zafra2023long}. Among these methods, \ac{FCS-MPC} has especially gained prominence due to its intuitive design, straightforward implementation, and high dynamic performance.\looseness=-1

Since the switch positions are discrete variables, \ac{FCS-MPC} results in an integer program. Its first variants utilized a horizon of one step and an exhaustive enumeration as solution approach~\cite{muller2003new,rodriguez2004predictive}. 
Subsequent studies revealed that long horizons play a crucial role in, e.g., mitigating harmonic distortion, and increasing the closed-loop stability margin~\cite{karamanakos_2014,geyer2014performance,geyer2014benefit}. Nonetheless, they render the exhaustive enumeration approach computationally infeasible. 
Whenever a linear prediction model is available, an efficient variant of a branch-and-bound algorithm called the sphere decoder has been shown to enable the use of long horizons within short sampling intervals~\cite{geyer_2014, geyer_2016}. \looseness=-1

Our interest lies in \ac{c1}, in which the torque and stator flux magnitude are controlled along their references.
\ac{c1} constitutes an important class of long horizon \ac{FCS-MPC}.
In particular, the prediction models have nonlinear input-output relations. In case the state-update function is linear except a nonlinear operation on the input itself, e.g., the absolute value operator found in the linearized neutral-point dynamics of a \ac{3L-NPC} converter, the work in \cite{liegmann_2017} has shown that one can lift the input space by augmenting it with the so-called ``pseudo-inputs" and implement the sphere decoder in the lifted space. When dealing with torque, the output function involves a nonlinear mapping from the states, rendering this approach inapplicable.
On the other hand, the slack variables associated with certain nonlinear constraints, e.g., on the switching frequency, have recently been integrated into the sphere decoder, while still maintaining its computational efficiency~\cite{hartmann2023switching}. Nevertheless, this approach is applicable only if a computationally lightweight, non-trivial, and provable lower bound exists for the future cost to be incurred from the slack variables; and no similar readily computable bound exists for any term involving torque.\looseness=-1

For the aforementioned reasons, the long horizon \ac{FCS-MPC} literature commonly adopts the concept of field orientation~\cite{hasse1968dynamischen}, see~\cite{liegmann_2017,stellato_2017,acuna2019impact,liegmann2021real} for example.
This notion enables the mapping of torque and flux references to stator current references through outer loops in the synchronous $\mathrm{dq}$-frame aligned with the angular position of the rotor flux vector~\cite{laczynski2009predictive,vargas2007predictive}. By doing so, this approach circumvents the nonlinear relation by capitalizing on the linearity of the current dynamics. Nevertheless, it is widely-recognized that controlling the torque directly maintains a consistent high-bandwidth torque control across various operating conditions~\cite{takahashi1986new,casadei2002foc,buja2004direct,miranda2009predictive}, especially when high switching frequencies are not viable, as in medium-voltage applications~\cite{geyer2008model,geyer2018algebraic}.\looseness=-1

An alternative approach to solve \ac{c1} is to develop a conventional branch-and-bound algorithm\cite{lawler1966branch,geyer2011computationally}. During transients, when a meaningful initial guess is unavailable, such heuristics may necessitate traversing a large portion of the decision tree before identifying the optimal solution. Since long horizons---thus, large trees---are desirable, this procedure needs to be stopped early to limit the computation time. One such stopping condition occurs when the number of nodes traversed exceeds a certain threshold. Even in the cases where the computationally efficient sphere decoder is applicable, similar limits are imposed~\cite{karamanakos2015suboptimal}. When the limit is reached, suboptimal early-stopping solutions would occasionally be applied. This could compromise control performance and even stability.\looseness=-1

In this paper, we will formulate the convex \ac{SDP} relaxation of \ac{c1}. 
Formulating \ac{SDP} relaxations of such polynomial optimization problems (with or without integer variables) has been extensively studied in the literature, e.g., consider the moment/sum-of-squares hierarchies of~\cite{lasserre2001global,parrilo2003semidefinite}. Furthermore, these relaxations have recently found practical applications in power electronics for offline computation of optimized pulse patterns~\cite{wachter2021convex}. 
An input sequence candidate will be extracted from the solution of the resulting~\ac{SDP} relaxation. While this sequence may only be approximately optimal for minimizing steady-state tracking errors, it often offers a superior solution during transients compared to the early-stopping branch-and-bound approach, and it can be obtained efficiently. Therefore, we will solve the SDP relaxation in parallel with the branch-and-bound algorithm and select the better sequence out of the two.\looseness=-1

Our contributions are as follows: 
\begin{itemize}
\item[\textit{(i)}] We propose a branch-and-bound algorithm to solve the long horizon \ac{FCS-MPC} problems with nonlinear system models as in torque and stator flux magnitude.  \looseness=-1

\item[\textit{(ii)}] We formulate the \ac{SDP} relaxation for \ac{c1}. To the best of our knowledge, this is the first work to do so. Moreover, the resulting formulation slightly differs from the existing works from the mathematical programming community on the \ac{SDP} relaxations of polynomial optimization, for example those over the Boolean hypercube, because our decision variables lie in the ternary hypercube~\cite{laurent2003comparison}. \looseness=-1

\item[\textit{(iii)}] We demonstrate in case studies that the extracted input sequence candidate can be advantageous during transients, especially when the node-limited branch-and-bound algorithm fails to find a meaningful input sequence.\looseness=-1
\end{itemize}

The paper is organized as follows. Section~\ref{sec:prel} presents the preliminaries, including the modeling, the \ac{c1} formulation, and the conventional branch-and-bound algorithm. Section~\ref{sec:relax} derives the~\ac{SDP} relaxation, and presents the proposed solution method. Finally, Section~\ref{sec:CaseStudies} showcases the performance of the method through numerical case studies.\looseness=-1

\textit{Notation:} Vectors are denoted with bold-face, e.g., $\b{x}\in\R^n$. Another vector could be formed by restricting the indices, e.g., $\b{x}_{i:j}=[{x}_{i}\,{x}_{i+1}\cdots{x}_{j}]^\top\in\R^{j-i}.$ Matrices are denoted capitalized, e.g., $\b{P}\in\R^{m\times n}$. Similar to the vector notation, $\b{P}_{i:j,k}$ restricts the matrix to its $k^\text{th}$ column and to the row range $i$~to~$j$. We denote the trace of $\b{P}\in \R^{n\times n}$ by $\trace({\b{P}})$, whereas $\diag({\b{P}})\in\R^n$ denotes the diagonal. The set of $n\times n$ real symmetric matrices is denoted by $\SM^n$.  The matrix $\b{P}\in \SM^n$
is \ac{PSD}, if $x^\top \b{P} x\geq 0$ for all $\b{x}\in\R^n$, and denoted by $\b{P}\succeq 0,$ i.e., $\b{P}\in \SM^n_+$. The set $\SM^n_+$ defines the \ac{PSD} cone. Finally, $\times$ is the cross-product, whereas $\odot$ is the Hadamard product.

\section{Preliminaries}\label{sec:prel}

\begin{figure}[t]
    \centering
\includegraphics[width=0.4875\textwidth]{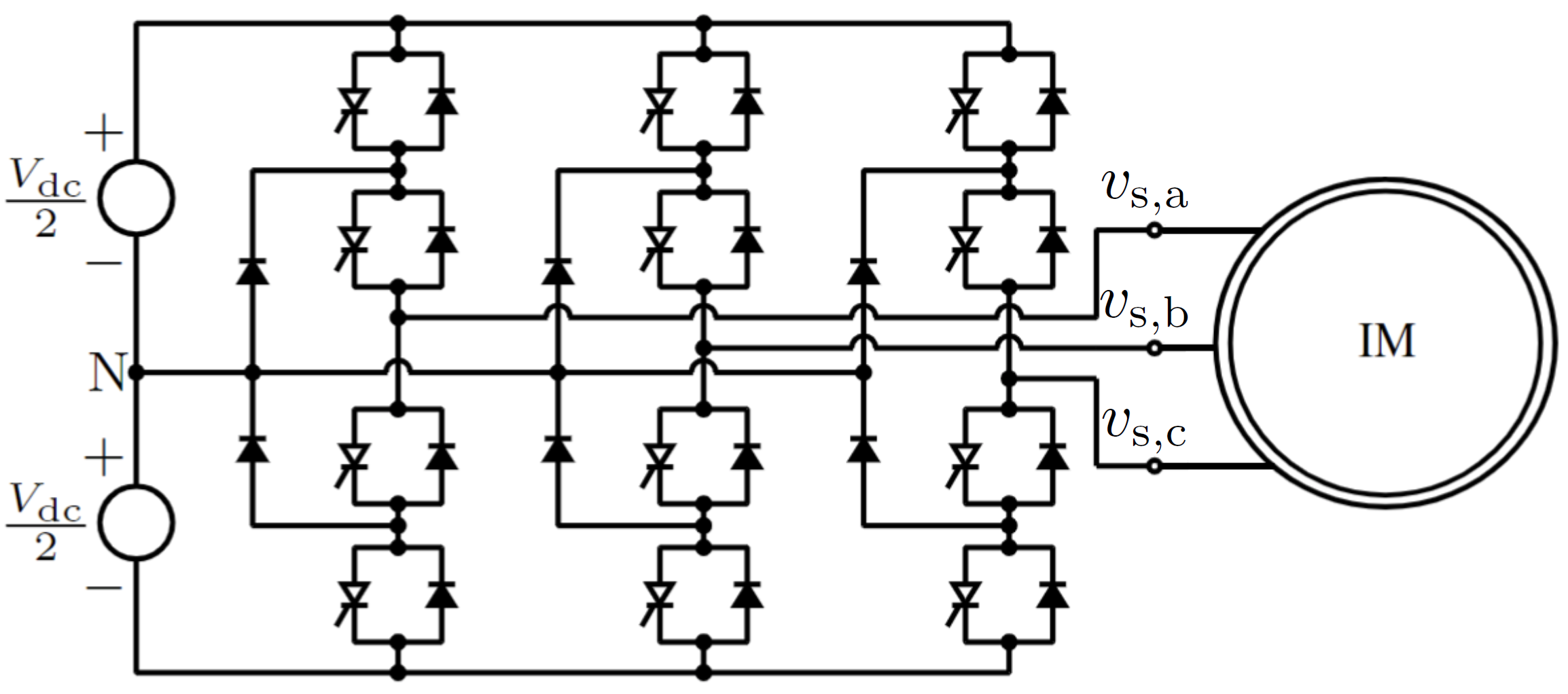}	
    \caption{NPC inverter with an induction motor and a fixed neutral-point potential.}
    \label{fig:3L-NPC_scheme}
\end{figure}
Consider a \ac{3L-NPC} converter connected to a squirrel cage \ac{IM}, as depicted in Figure~\ref{fig:3L-NPC_scheme}. The half dc-link voltages $\frac{V_\dc}{2}$ are realized by ideal voltage sources, thus fixing the neutral-point potential to zero. The reader is referred to~\cite{liegmann_2017, Jin_2022} for more details on neutral-point balancing. \looseness=-1

Throughout the paper, the dynamical equations are based on the $\alpha\beta$-reference frame (also known as the {stationary/orthogonal frame}).
It is denoted by $\b{\xi}_{\alpha\beta} = [\xi_\alpha~\xi_\beta]^\top$, where
$\b{\xi}_{\alpha\beta} = \b{K}\b{\xi}_{abc}$
is the {Clarke} transformation with $\b{K} = \frac{2}{3}{\tiny\begin{bmatrix}
			1 & -\frac{1}{2} & -\frac{1}{2} \\
			0 & \frac{\sqrt{3}}{2} & -\frac{\sqrt{3}}{2}
		\end{bmatrix}}$.
Similarly, the inverse transformation is defined as
$\b{\xi}_{abc}  = \b{K}^{-1}\b{\xi}_{\alpha\beta},$
where the (pseudo) inverse of $\b{K}$ is
$\b{K}^{-1} = \frac{3}{2}\b{K}^{\top}.$

It is convenient to model the states of \ac{IM} in terms of its stator and rotor flux linkage vectors $\b{\psi}_{\s, \alpha\beta} = [\psi_{\s, \alpha}\ \psi_{\s, \beta}]^\top$ and $\b{\psi}_{\r, \alpha\beta} = [\psi_{\r, \alpha}\ \psi_{\r, \beta}]^\top$, respectively. The switch positions are denoted by $\b{u}_\abc=[{u}_\A\ {u}_\B\ {u}_\C]^\top\in \{-1,0,1\}^3$. The stator voltages $\b{v}_{\s,\abc} = [v_{\s,\A}\ v_{\s,\B}\ v_{\s,\C}]^\top= \frac{V_\dc}{2}\b{u}_\abc \in \{-\frac{V_\dc}{2},\ 0,\ \frac{V_\dc}{2}\}^3$ are the output converter voltages at the ac-side. We drop the subscript $\abc$ from $\b{u}_\abc$ and simply write $\b{u}$.  \looseness=-1 


\vspace{0cm}
\subsection{Physical system model}
\vspace{0cm}
The continuous-time flux linkage dynamics of \ac{IM} are: 
{\medmuskip=0.001mu\thinmuskip=0.001mu\begin{equation*}
    \begin{split}
    \frac{d\b{\psi}_{\s,\alpha\beta}(t)}{dt} &= -R_\s \frac{X_\r}{D} \b{\psi}_{\s,\alpha\beta}(t) + R_\s \frac{X_\m}{D} \b{\psi}_{\r,\alpha\beta}(t) + \frac{V_\dc}{2}\b{K}\b{u}(t), \\
    \frac{d\b{\psi}_{\r,\alpha\beta}(t)}{dt} &= R_\r \frac{X_\m}{D} \b{\psi}_{\s,\alpha\beta}(t)- \Big(R_\r \frac{X_\s}{D} - \omega_\r 
    \begin{bmatrix}
    0 & -1 \\
    1 & 0 
    \end{bmatrix}\Big)
    \b{\psi}_{\r,\alpha\beta}(t), 
    \end{split}
\end{equation*}}where $R_\s$, $R_\r$ are the stator and rotor resistances, respectively; $X_\s$, $X_\r$, $X_\m$ are the stator, rotor, and main inductances, respectively; and $D = X_\s X_\r - X_\m^2$~\cite{krause2002analysis}, \cite[\S 4.3]{book_geyer}. For the sake of simplicity, via time-scale separation, we treat the rotor's angular speed $\omega_\r$ as a parameter, and discussions on the outer speed controller are considered out-of-scope for this paper. \looseness=-1

Using exact Euler discretization with the control sampling interval $T_\C$, the discrete-time model becomes
{\medmuskip=0.001mu\thinmuskip=0.001mu\begin{equation}\label{eq:disc_dynamics}
    \begin{split}
    \b{\psi}_{\s,\alpha\beta}(k+1) &= \b{A}_{11} \b{\psi}_{\s,\alpha\beta}(k) + \b{A}_{12} \b{\psi}_{\r,\alpha\beta}(k) + \b{B}_1 \b{u}(k), \\
    \b{\psi}_{\r,\alpha\beta}(k+1) &= \b{A}_{21} \b{\psi}_{\s,\alpha\beta}(k) + \b{A}_{22} \b{\psi}_{\r,\alpha\beta}(k),
    \end{split}
\end{equation}}where $k\in \mathbb{N}$ is the discrete-time index, and $\b{A}_{11}$, $\b{A}_{12}$, $\b{A}_{21}$, $\b{A}_{22}$, and $\b{B}_1$ can be inferred from the continuous dynamics.\looseness=-1

The control task is to track the references of the electromagnetic torque $T_\e$ and the stator flux magnitude $\Psi_\s$. The torque is\looseness=-1
{\medmuskip=0.01mu\thinmuskip=0.01mu\begin{align}
    T_\e(k) &= \frac{1}{\mathrm{pf}} \frac{X_\m}{D} \b{\psi}_{\r,\alpha\beta}(k) \times \b{\psi}_{\s,\alpha\beta}(k)\nonumber \\
&= T_\factor \b{\psi}_{\r,\alpha\beta}(k)^\top \b{J} \b{\psi}_{\s,\alpha\beta}(k). \label{eq:T_e}
\end{align}}The power factor is $\mathrm{pf}=P_\mathrm{rat}/S_\mathrm{rat}$, where $P_\mathrm{rat}$ and $S_\mathrm{rat}$ are the rated real and apparent power, respectively.
Moreover, $T_\factor = \frac{1}{\mathrm{pf}} \frac{X_\m}{D}$, and $ \b{J} = {\tiny\begin{bmatrix}
    0 & 1 \\ -1 & 0
\end{bmatrix}}$. The stator flux magnitude is
\begin{equation}\label{eq:Psi_s}
    \Psi_\s(k) = ||\b{\psi}_{\s,\alpha\beta}(k)||_2 = \sqrt{\psi_{\s,\alpha}(k)^2 + \psi_{\s,\beta}(k)^2}.
\end{equation}

\vspace{0cm}
\subsection{Finite-control-set model predictive control problem}
\vspace{0cm}
\begin{figure}[t]
    \includegraphics[width=0.475\textwidth]{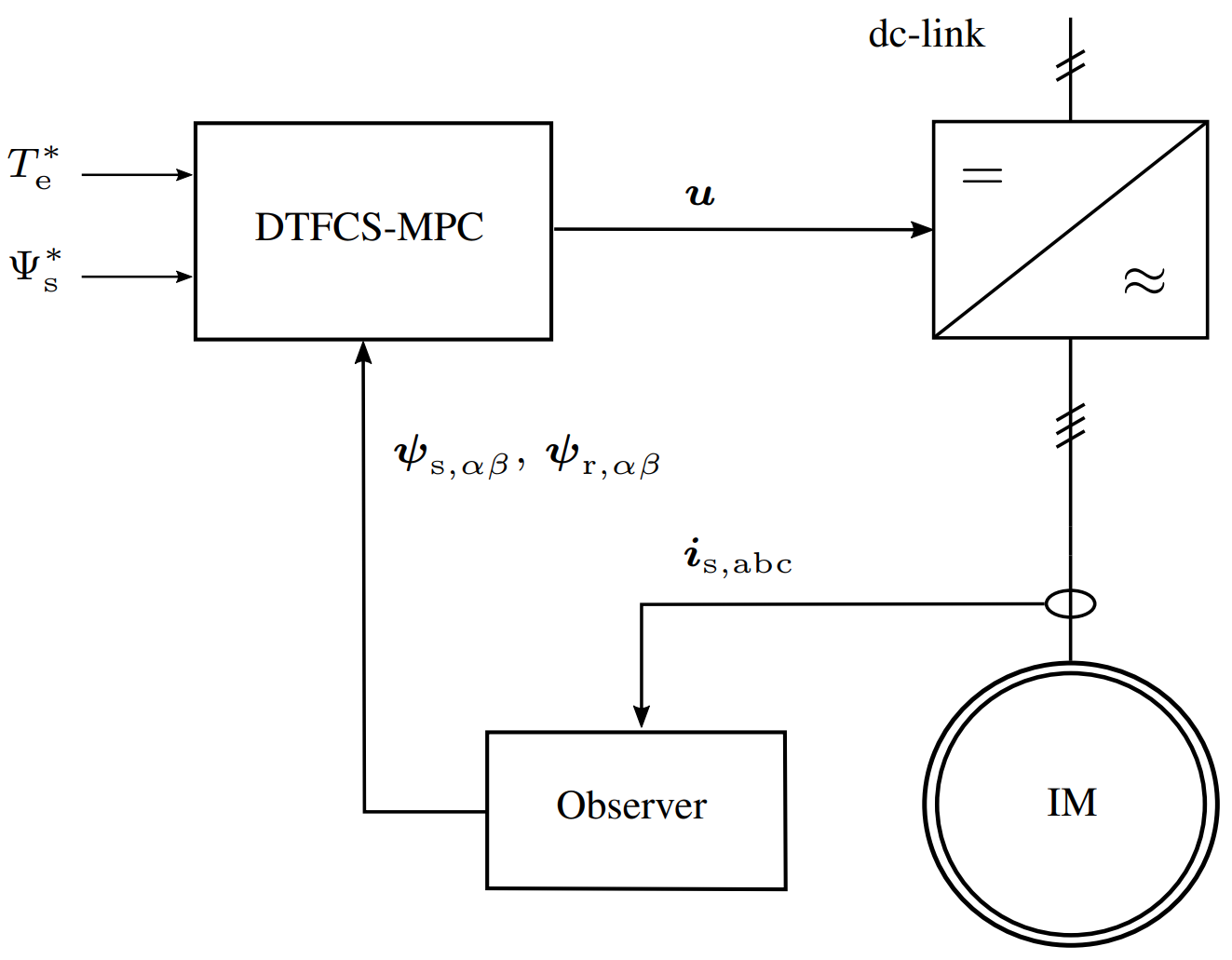}	
    \caption{A high-level \ac{c1} block diagram.}
    \label{fig:MPC_block-diagram}
\end{figure} 

A high-level block diagram of \ac{c1} is shown in Figure~\ref{fig:MPC_block-diagram}. The prediction horizon is denoted by $N \in \mathbb{N}$. It tracks both the torque and the stator flux magnitude references, usually generated by certain outer loops, e.g., the speed controller. A study on the impact of the observer design is considered out-of-scope for this paper; kindly refer to~\cite[\S 5.3]{sul2011control} for more details.\looseness=-1

Define a compact dynamical system with states  $\b{x}(k) = [\b{\psi}_{\s,\alpha\beta}(k)^\top\ \b{\psi}_{\r,\alpha\beta}(k)^\top]^\top$, input $\b{u}(k)$, and output $\b{y}(k) = [T_\e(k)\ \Psi_\s^2(k)]^\top$.\footnote{Squaring the magnitude in~\eqref{eq:Psi_s} turns it into a polynomial.} The state dynamics are given by $\b{x}(k+1) = \b{A}\b{x}(k) + \b{B}\b{u}(k),$ where $\b{A}$ and $\b{B}$ can be inferred from~\eqref{eq:disc_dynamics}. Note that the state-update function is linear. The output function is nonlinear, and it follows from the definitions of $T_\e(k)$ and $\Psi_\s(k)$ in~\eqref{eq:T_e} and~\eqref{eq:Psi_s}, respectively. \looseness=-1

Let $e_\T(k) = T_\e^* -  y_1(k)$ and $e_\Psi(k) = \Psi^{*2}_\s -  y_2(k)$ be the tracking errors for torque and stator flux magnitude, respectively.\footnote{For the sake of simplicity, the references are time-invariant.} The objective of \ac{c1} is \looseness=-1
{\medmuskip=0.01mu\thinmuskip=0.01mu\begin{align*}
    f_N(\b{U}(k)) = &\sum_{\ell=1}^{N_\p} \lambda_\T\left( e_\T(k+\ell) \right)^2 \\
    &\quad+ (1 - \lambda_\T) \left( e_\Psi(k+\ell) \right)^2 + \lambda_\U ||\Delta \b{u}(k+\ell-1)||_2^2,
\end{align*}}where $\ell\in\{1,\ldots,N\}$ is used to iterate over the prediction horizon. The parameter $\lambda_\T$ sets the trade-off between the two objectives. Additionally, a penalty with the weight $\lambda_\U$ on the switching transitions $\Delta \b{u}(k) = \b{u}(k) - \b{u}(k-1)$ is included to limit the switching frequency. Not penalizing control effort results in a deadbeat behaviour~\cite[\S IV.A]{karamanakos2019guidelines}. Here, we also introduce the full-horizon switching sequence $\b{U}(k) = [\b{u}(k)^\top\ \b{u}(k+1)^\top\ \cdots\ \b{u}(k+N-1)^\top]^\top$. \looseness=-1

\ac{c1} solves the following optimization problem at time step $k$:
{\medmuskip=0.1mu\thinmuskip=0.1mu \begin{equation}\label{eq:c1-opt} \tag{$\mathcal{P}$}
   \begin{split}
    f^\opt_N=&\min_{\b{U}(k)}\ f_N(\b{U}(k)) \\
    &\mathrm{s.t.}\ \b{x}(k+\ell) = \b{A}\b{x}(k+\ell-1) + \b{B}\b{u}(k+\ell-1), \\
     & \quad\ \ y_1(k+\ell) = T_\factor \b{x}_{3:4}(k+\ell)^\top \b{J} \b{x}_{1:2}(k+\ell), \\
    & \quad\ \ y_2(k+\ell) = ||\b{x}_{1:2}(k+\ell)||_2^2, \\
    & \quad\ \ \b{u}(k+\ell-1) \in \{-1,\ 0,\ 1\}^{3},\quad \forall \ell. 
    \end{split}
\end{equation}}After solving~\eqref{eq:c1-opt} for its optimal solution $\b{U}^\opt(k)$, \ac{c1} applies the first input $\b{u}^\opt(k)$ in a receding horizon fashion. As will be discussed further in Section~\ref{sec:relax}, \eqref{eq:c1-opt} is a polynomial optimization problem over integer variables. The objective $ f_N(\cdot)$ is of degree~4. This class is known to be NP-hard~\cite{hastad1996clique}. Let $\mathcal{U}_N = \{-1,\ 0,\ 1\}^{3N}$ denote all possible input sequences for~\eqref{eq:c1-opt}. Since $\mathrm{card}(\mathcal{U}_N) = 27^N$, enumeration is computationally infeasible for long horizons, i.e., large $N$.
\looseness=-1

\subsection{A branch-and-bound algorithm for \ac{c1}} \label{sec:branch-and-bound}

We briefly recall a conventional branch-and-bound algorithm.
Let $\b{U}^{\ini}(k)$ denote the initial solution with the objective value $f^\ini_N = f_N(\b{U}^\ini(k))$. 
For example, it can be obtained by an educated guess using a shifted version of the optimal solution of the previous time step $k-1$, i.e., $\b{U}^\opt(k-1)$.
Starting with this being the incumbent solution, the algorithm then systematically searches for the optimal solution. Whenever a better input sequence is identified, it updates its incumbent solution. During its exploration, it can shrink the search space based on the objective value of the incumbent solution.\looseness=-1

This exploration procedure with shrinking relies on the additive nonnegative structure of the objective. Instead of directly calculating $f_N(\b{U}(k))$ for a full candidate solution, we create a solution tree of depth $N$, where at each node, we select the input $\b{u}(k+\ell-1)$ of the time step $k+\ell-1$ at the depth level $\ell\in\{1,\ldots,N\}$. As we choose new inputs at later time steps, we can iteratively increment the objective to obtain \looseness=-1
{\medmuskip=0.01mu\thinmuskip=0.01mu\begin{equation} \label{eq:incremental_cost}
\begin{split}
    f_\ell = &\sum_{j=1}^{\ell} \lambda_\T\left(e_\T(k+j) \right)^2\\&\quad+ (1 - \lambda_\T) \left(e_\Psi(k+j) \right)^2 + \lambda_\U ||\Delta \b{u}(k+j-1) ||_2^2.
    \end{split}
\end{equation}}Whenever the incremental objective $f_\ell$ exceeds the objective value of the incumbent solution, the sub-tree can simply be discarded. \looseness=-1

With this in mind, the branch-and-bound algorithm is presented in Algorithm~\ref{alg:1}. To initialize, we let $f_0 = 0$, $\ell=1$, $\b{U} = \b{U}^\inc = \b{U}^\ini(k)$, $f_N^\inc = f^\ini_N$, $\b{u}^{\prev} = \b{u}(k-1)$, $\b{x}^{\prev} = \b{x}(k-1)$, $n_\p = 0$. After deciding upon $\b{u}(k+\ell-1)$, the resulting incremental objective~\eqref{eq:incremental_cost} is compared to the best so far, i.e., the incumbent solution, in Line~\ref{line:incumbent_solution}. We continue to traverse if $f_\ell < f_N^\inc$ holds. Otherwise, the sub-tree is pruned. This is repeated until $f_N < f_N^\inc$ is found at a leaf node, and then the incumbent solution is updated. The algorithm starts back-tracking. The optimal solution $\b{U}^\opt(k)$ is found and its optimality is verified if all the branches have been cut. Otherwise, the algorithm stops early whenever the node limit, $n_{\p,\max}$, is hit.\looseness=-1

The branch-and-bound algorithm has a worst-case exponential time performance. Nevertheless, it is efficient during steady-state operation. During transients, however, the initial solution may not be meaningful. Moreover, in contrast to the sphere decoder of~\cite{geyer_2014}, which was proposed for linear systems with integer manipulated variables, it does not incorporate the future cost incurred by an input decision at an earlier time step. In Algorithm~\ref{alg:1}, the cost is simply incremented by the stage cost incurred at this time step. Thus, it might need to traverse a large portion of the tree to find the optimal solution and to verify its optimality. Thus, it is common practice to impose a time limit, e.g., a node limit as in Line~\ref{line:node_limit}, which might then result in the algorithm returning suboptimal solutions.

\begin{algorithm}[t]
\caption{Conventional branch-and-bound algorithm:} 
\label{alg:1}
\begin{algorithmic}[1]
    \State \textbf{Function:} branch\_and\_bound
    \State \textbf{Input:} $f_{\ell-1}$, $f_N^\inc$, $\ell$, $\b{U}$, $\b{U}^\inc$, $\b{x}_\prev$, $\b{u}_\prev$, $n_\p$
    \State \textbf{Output:} $\b{U}^\inc$
    \If{$n_\p > n_{\p,\max}$}\label{line:node_limit}\ \textbf{return} $\b{U}^\inc$
    \Else\ $n_\p \gets n_\p + 1$
    \EndIf 
    \ForAll{$\b{u}\in\mathcal{U}_1$}
        \State $\b{x}(k+\ell) \gets \b{A}\b{x}_\prev + \b{B}\b{u}$
        \State $y_1(k+\ell) \gets T_\factor \b{x}_{3:4}(k+\ell)^\top \b{J} \b{x}_{1:2}(k+\ell)$
        \State $y_2(k+\ell) \gets ||\b{x}_{1:2}(k+\ell)||_2^2$
        \State $f_\ell \gets f_{\ell-1} + \lambda_\T(e_\T(k+\ell))^2 + (1 - \lambda_\T) (e_\Psi(k+\ell))^2 + \lambda_\U || \b{u} - \b{u}_\prev ||_2^2$
        \If{$f_\ell < f_N^\inc$} \label{line:incumbent_solution}
            \State $\b{U}_{3\ell-2:3\ell} \gets \b{u}$
            \If{$\ell < N$}
                \State $\b{U}^\inc\gets$ branch\_and\_bound($f_\ell$, $f_N^\inc$, $\ell+1$, $\b{U}$, $\b{U}^\inc$, $\b{x}(k+\ell)$, $\b{u}$, $n_\p$)
            \Else
                \State $\b{U}^\inc \gets \b{U}$
                \State $f_N^\inc \gets f_\ell$
            \EndIf
        \EndIf
    \EndFor 
    \State \textbf{return} $\b{U}^\inc$
\end{algorithmic}
\end{algorithm}

\section{Exploiting \ac{SDP} relaxations to solve \ac{c1}}\label{sec:relax}
This section formulates the \ac{SDP} relaxation of~\eqref{eq:c1-opt}. Different from~\eqref{eq:c1-opt}, \ac{SDP}s are well-studied conic convex optimization problems with many solvers developed to address them.\footnote{In terms of complexity, \ac{SDP}s are considered to be in PSPACE~\cite{tarasov2008semidefinite,o2017sos}.} To this end, we first formulate $f_N(\b{U}(k))$ explicitly as a function of $\b{U}(k)$, the initial state $\b{x}(k)$ and the previous input $\b{u}(k-1)$. \looseness=-1
From~\eqref{eq:disc_dynamics}, derive the maps
\begin{equation*}
    \begin{split}
       \b{x}_{1:2}(k+\ell) & = \b{{\psi}}_\s(k+\ell) = \b{\Gamma}_{\s,\ell} \b{x}(k) + \b{\Upsilon}_{\s,\ell} \b{U}(k),\\
    \b{x}_{3:4}(k+\ell) & =\b{{\psi}}_\r(k+\ell) = \b{\Gamma}_{\r,\ell} \b{x}(k) + \b{\Upsilon}_{\r,\ell} \b{U}(k),
    \end{split}
\end{equation*}
where $\b{\Gamma}_{\s,\ell}$, $\b{\Gamma}_{\r,\ell}$, $\b{\Upsilon}_{\s,\ell}$, $\b{\Upsilon}_{\r,\ell}$ are relegated to Appendix~\ref{app:matrices}. Notice they are $k$-invariant, and can thus be computed offline. \looseness=-1

Now, augment the full-horizon variable with $1$ as follows $\b{\tilde{U}}(k) = [1\ \b{U}(k)^\top]^\top$. Given the definitions of $T_\e(k)$ and $\Psi_\s(k)$ in~\eqref{eq:T_e} and~\eqref{eq:Psi_s}, respectively, we have
\begin{equation}  \label{eq:pol}
    \begin{split}
    y_1(k+\ell) &= T_\e(k+\ell) = \b{\tilde{U}}(k)^\top \b{Q}_\ell(k) \b{\tilde{U}}(k),\\
    y_2(k+\ell) &= \Psi_\s^2(k+\ell) = \b{\tilde{U}}(k)^\top \b{W}_\ell(k) \b{\tilde{U}}(k),
    \end{split}
\end{equation}where $\b{W}_\ell(k)$ and $\b{Q}_\ell(k)$ are functions of $\b{x}(k)$ and relegated to Appendix~\ref{app:matrices}.\looseness=-1 

For switching transitions, we have \looseness=-1
\begin{equation}\label{eq:pol_du}
	\Delta u_z(k+\ell-1) =
	\b{Z}_{z,\ell}(k)
    \b{\tilde{U}}(k),\ \forall z\in\{\A,\B,\C\},
\end{equation}
where $\b{Z}_{z,\ell}$ is relegated to Appendix~\ref{app:matrices}. Its first instance $\b{Z}_{z,1}$ is a function of $\b{u}(k-1)$, and the rest can be computed offline. \looseness=-1

Using the definitions above, the objective of~\eqref{eq:c1-opt} becomes
{\medmuskip=0.01mu\thinmuskip=0.01mu\begin{equation}\label{eq:c1obj}
\begin{split}
    f_N(&\b{{U}}(k)) = \lambda_T \sum_{\ell = 1}^{N} \left(T^*_\e-\b{\tilde{U}}(k)^\top \b{Q}_\ell(k) \b{\tilde{U}}(k) \right)^2 \\
	&\ + (1-\lambda_T) \sum_{\ell = 1}^{N} \left(\Psi_s^{*2}-\b{\tilde{U}}(k)^\top \b{W}_\ell(k) \b{\tilde{U}}(k)\right)^2  \\
	&\ + \lambda_u \sum_{\ell = 1}^{N}\ \sum_{z \in \{\A,\B,\C\}} \b{\tilde{U}}(k)^\top \b{Z}_{z,\ell}(k)^\top\b{Z}_{z,\ell}(k) \b{\tilde{U}}(k).
 \end{split}
\end{equation}}With this, we restate~\eqref{eq:c1-opt} without the state variables
{\medmuskip=0.01mu\thinmuskip=0.01mu\begin{equation} \label{eq:c1_opt_}
\begin{split}
	f^\opt_N=&\min_{\b{{U}}(k)}\ f_N(\b{{U}}(k)) \\
	&\begin{array}{llr}
		\mathrm{s.t.} & \b{U}(k) \in \mathcal{U}_N. \\
	\end{array}
\end{split}
\end{equation}}It is now straightforward to see that we have been working with a degree-4 polynomial over a ternary hypercube.

\subsection{SDP relaxation of \ac{c1} optimization problem}

With a change of variables, and the identity $\b{x}^\top\b{P}\b{x}=\trace(\b{P}\b{x}\b{x}^\top)$, \eqref{eq:c1obj} is equivalent to
{\medmuskip=0.01mu\thinmuskip=0.01mu\begin{equation}\label{eq:newobj}
\begin{split}
 f_N^{\b{\Theta}} (&\b{\Theta}(k)) = \lambda_T \sum_{\ell = 1}^{N} \left(T^*_\e-\trace\left( \b{Q}_\ell(k)\b{\Theta}(k)\right) \right)^2 \\
	&\quad + (1-\lambda_T) \sum_{\ell = 1}^{N} \left(\Psi_s^{*2}-\trace\left(  \b{W}_\ell(k)\b{\Theta}(k)\right) \right)^2\\
	&\quad + \lambda_u \sum_{\ell = 1}^{N}\sum_{z \in \{\A,\B,\C\}} \trace\left(  \b{Z}_{z,\ell}(k)^\top\b{Z}_{z,\ell}(k)\b{\Theta}(k)\right),
 \end{split}
\end{equation}}where $\b{\Theta}(k) = \b{\tilde{U}}(k)\b{\tilde{U}}(k)^\top\in \{-1,\ 0,\,1\}^{3N+1\times3N+1}$. 
The variable $\b{\Theta}(k)$ lies in the lifted space and the relation $\b{\Theta}(k) = \b{\tilde{U}}(k)\b{\tilde{U}}(k)^\top$ captures its dependence to the original variable, so that $f_N^\Theta(\b{\Theta}(k)) = f_N(\b{U}(k))$.
With \eqref{eq:newobj}, \eqref{eq:c1_opt_} can be reformulated as the following
{\medmuskip=0.01mu\thinmuskip=0.01mu\begin{equation} \label{eq:c1_opt_Theta}
\begin{split}
	f^\opt_N=&\min_{\b{\Theta}(k)}\ f_N^{\b{\Theta}} (\b{\Theta}(k)) \\
	&\begin{array}{llr}
		\mathrm{s.t.} & \b{\Theta}(k) =  \b{\tilde{U}}(k)\b{\tilde{U}}(k)^\top,\\ & \b{U}(k) \in \mathcal{U}_N. \\
	\end{array}
\end{split}
\end{equation}}The constraint $ \b{\Theta}(k) =  \b{\tilde{U}}(k)\b{\tilde{U}}(k)^\top$ is equivalent to having that $\b{\Theta}(k)$ is \ac{PSD} (i.e., has a Cholesky decomposition), and is of rank 1. Reflecting on this observation, we revisit a relaxation method, first proposed and analyzed for the Max-Cut problem by the seminal works of~\cite{delorme1993laplacian,goemans1995improved}. 
This method substitutes the constraints of type $\b{\Theta}(k) =  \b{\tilde{U}}(k)\b{\tilde{U}}(k)^\top$ with the \ac{PSD} cone and a rank constraint. Subsequently, we can disregard the rank and integer constraints to obtain an \ac{SDP} relaxation.  

First, we obtain the equivalent reformulation:
\begin{equation} \label{eq:c1_opt_rank_constraint}
\begin{split}
	f^\opt_N=&\min_{\b{\Theta}(k)}\  f_N^{\b{\Theta}} (\b{\Theta}(k)) \\
	&\begin{array}{llr}
		\mathrm{s.t.} & \b{\Theta}(k) \succeq 0,\, \Theta_{1,1}(k) = 1,\\
        &\b{\Theta}(k) \in \{-1,\ 0,\,1\}^{3N+1\times3N+1},\\
        &\mathrm{rank}(\b{\Theta}(k)) = 1, \\
	\end{array} 
 \end{split}
\end{equation}
where $\Theta_{1,1}(k) = 1$ imposes the augmented 1.\footnote{$\mathrm{diag}(\b{\Theta}(k)) \geq 0$ holds by definition for the \ac{PSD} cone.} 
We now drop the nonconvex rank and integer constraints as follows: 
\begin{equation} \label{eq:c1_opt_relax}\tag{$\mathcal{P}_\mathrm{sdp}$}
\begin{split}
	f^\opt_N\geq &\min_{\b{\Theta}(k)}\  f_N^{\b{\Theta}} (\b{\Theta}(k)) \\
	&\begin{array}{llr}
		\mathrm{s.t.} & \b{\Theta}(k) \succeq 0,\, \Theta_{1,1}(k) = 1,\\
        &\b{\Theta}(k) \in [-1,\,1]^{3N+1\times3N+1}. \\
	\end{array} 
 \end{split}
\end{equation}
Let $\b{\Theta}^\rel(k)$ denote its optimal solution.
As a remark, the problem above is not in the standard form of an \ac{SDP}, as it involves a convex quadratic function of a trace of a matrix in~\eqref{eq:newobj}. Since trace is a linear map, this problem remains convex, and it can be reformulated with a linear objective using an epigraph representation and the Lorentz cone.\footnote{The Lorentz cone is linearly isomorphic to the $2\times2$ \ac{PSD} cone~\cite{ben2001lectures}, and can be handled by \ac{SDP} solvers.}
\ac{SDP}s are convex optimization problems that can be handled by a variety of off-the-shelf solvers. For case studies, we will use a first-order solver, SCS~\cite{o2016conic}, to obtain a solution fast without stringent precision requirements. Some of the other solvers that can handle \ac{SDP}s can be listed as~\cite{sturm1999,goulart2024clarabel,mosek2024}.
\looseness=-1

After solving~\eqref{eq:c1_opt_relax}, and obtaining $\b{\Theta}^\rel(k)$ in the lifted space, we need to extract an input sequence candidate. 
Towards this, we provide the following four options:
\begin{enumerate}
	\item \textit{First column:} We take the first column of $\b{\Theta}^\rel(k)$ and round it, i.e., 
		$\b{U}^{\mathrm{sdp}}(k) = \mathrm{round}\left(\b{\Theta}^\rel_{2:3N,1}(k)\right).$
	\item \textit{Diagonal:} We take the square-root of the diagonal, multiply it with the sign of the first column, and round the result, i.e., 
		$\b{U}^{\mathrm{sdp}}(k) = \mathrm{round}\left(\mathrm{sign}(\b{\Theta}^\rel_{2:3N,1}(k))\odot\sqrt{\mathrm{diag}(\b{\Theta}^\rel_{2:3N,2:3N}(k))}\right).$
	\item \textit{Maximal Eigenvector:} From the \ac{ED} of $\b{\Theta}^\rel(k)$ we round the eigenvector with the largest eigenvalue, i.e., 
		$\b{U}^{\mathrm{sdp}}(k) = \mathrm{round}\left(\sqrt{\lambda_1}\b{v}_1\right),$
	 where $\lambda_1$ is the largest eigenvalue of $\b{\Theta}^\rel(k)$ and $\b{v}_1$ the corresponding eigenvector.
        \item \textit{Weighted Eigenvectors:} From the \ac{ED} of $\b{\Theta}^\rel(k)$ we round the eigenvalue-weighted sum of the eigenvectors, i.e., 
		$\b{U}^{\mathrm{sdp}}(k) = \mathrm{round}\left(\sum_j\sqrt{\lambda_j}\b{v}_j\right).$
\end{enumerate}
Our case studies will rely on the computationally-simplest first option, as no significant control performance differences have been observed during transients.\footnote{The \ac{ED} has the complexity $\mathcal{O}(N^3)$. In practice, it can be computed with very good precision via several iterations of the Power Method, each with the complexity $\mathcal{O}(N^2)$.}
Note that for this class of problems, no theoretical approximation guarantee is known for the extracted candidates~\cite{henrion2005detecting}.
This remains open for future research. For quadratic optimization problems over the Boolean hypercube, for example, \cite{nesterov1998semidefinite} provides an approximation guarantee by generalizing the analysis of the GW algorithm of \cite{goemans1995improved}.  Nevertheless, the term $$f_N(\b{{U}}^{\mathrm{sdp}}(k)) -f_N^{\b{\Theta}}(\b{\Theta}^\rel(k)),$$ bounds the optimality gap, since $$f_N(\b{{U}}^{\mathrm{sdp}}(k))\geq f^\opt_N \geq f_N^{\b{\Theta}}(\b{\Theta}^\rel(k)).$$
Moreover, regardless of whether the situation is transient or steady-state, \ac{SDP}s are convex programs that can be solved within a reasonable and consistent time frame. 

Next, we present our proposed approach that exploits the advantages of both the \ac{SDP} and the branch-and-bound.
\looseness=-1

\subsection{Proposed algorithm}
We propose solving the \ac{SDP} relaxation in \eqref{eq:c1_opt_relax} in parallel to Algorithm~\ref{alg:1}. In the final step, we choose the best of the early-stopping solution from Algorithm~\ref{alg:1} and the sequence extracted from the \ac{SDP} solution in the lifted space. This procedure is formalized in Algorithm~\ref{alg:2}, where the \textit{first column} extraction method is adopted.\looseness=-1

This proposal is constructed to harness the benefits of both the \ac{SDP} and the branch-and-bound algorithm. It achieves optimality during the steady-state, given that the branch-and-bound concludes prior to hitting the node limit. The \ac{SDP} relaxation is also solved fast consistently in parallel and provides a meaningful close-to-optimal input sequence even during transient events.\looseness=-1

\begin{remark}
There exist alternative strategies for leveraging the \ac{SDP} relaxations that we are not examining in this study.
One can use the best of the sequence extracted from the \ac{SDP} and the educated guess to warm-start the branch-and-bound algorithm.
However, this would be computationally wasteful as the \ac{SDP} and the branch-and-bound would be run sequentially. 
Simultaneously solving \ac{SDP}s during the branch-and-bound iterations is yet another option, as the \ac{SDP} could provide a lower bound on the future branch cost. However, this might not be realistic from a computational standpoint. 
\end{remark}

\newlength\myindent
\setlength\myindent{2em}
\newcommand\bindent{%
  \begingroup
  \setlength{\itemindent}{\myindent}
  \addtolength{\algorithmicindent}{\myindent}
}
\newcommand\eindent{\endgroup}

\begin{algorithm}[t]
\caption{Branch-and-bound with the SDP relaxation}
\label{alg:2}
\begin{algorithmic}[1]
    \State \textbf{Function:} branch\_and\_bound\_with\_SDP
    \State \textbf{Input:} $\b{x}(k-1), \b{u}(k-1), \b{U}^\opt(k-1)$
    \State \textbf{Output:} $\b{U}^\inc, n_\p$ \\ 
    Start Parallel computation:
    \State \textbf{Thread 1} (Branch-and-bound)
    \State $\b{U}^\mathrm{ed} \gets [\b{U}_{4:3N}^\opt(k-1)^\top \b{U}_{3N-2:3N}^\opt(k-1)^\top]^\top$
    \State  $f^\mathrm{ed} \gets f_N(\b{{U}}^\mathrm{ed})$
    \State  $f_0 \gets 0,\ \ell \gets 1,\ n_\p \gets 0$
    \State $\b{U}^{\mathrm{b\&b}}\gets\text{ branch\_and\_bound}(f_0,f^\mathrm{ed},\ell,\b{U}^\mathrm{ed},$
    $\b{U}^\mathrm{ed},\b{x}(k-1),\b{u}(k-1),n_\p)$
    \State \textbf{End Thread 1}
    \State \textbf{Thread 2} (SDP relaxation)
    \State  Solve~\eqref{eq:c1_opt_relax} for $\b{\Theta}^\mathrm{sdp}(k)$
    \State  $\b{U}^\mathrm{sdp} \gets \mathrm{round}(\b{\Theta}^\mathrm{sdp}_{2:3N,1}(k))$
    \State \textbf{End Thread 2}
    \If{$f_N(\b{U}^{\mathrm{b\&b}}) < f_N(\b{U}^{\mathrm{sdp}})$}
    \State $\b{U}^\inc \gets \b{U}^{\mathrm{b\&b}}$
    \Else
    \State $\b{U}^\inc \gets \b{U}^{\mathrm{sdp}}$
    \EndIf
    \State \textbf{return} $\b{U}^\inc$
\end{algorithmic}
\end{algorithm}

\section{Case Studies}\label{sec:CaseStudies}
First, we present a case study illustrating the variability in the efficiency of Algorithm~\ref{alg:1}, due to the initial guess and the operational conditions, when no node limits are enforced. Subsequently, we demonstrate through another case study how Algorithm 1 could potentially fail to obtain a meaningful solution during transients. Algorithm~\ref{alg:2} then effectively remedies this issue. \looseness=-1

\subsection{Parameters and simulation setting}

\begin{table}[t]
	\centering
	\caption{Physical parameters.}
	\label{tab:sim_parameters}
	\begin{tabular}{cccc}
		\hline
		\multicolumn{2}{c}{Rated values} & \multicolumn{2}{c}{Machine parameters [pu]} \\
		\hline
		Voltage & $\SI{3300}{\volt}$ & Stator resistance & $R_s = 0.0108$ \\ 
		Current & $\SI{356}{\ampere}$ & Rotor resistance & $R_r = 0.0091$ \\ 
		Real power & $\SI{1.587}{\mega\watt}$ & Stator leak. react. & $X_{ls} = 0.1493$ \\ 
		Apparent power & $\SI{2}{\mega\voltampere}$ & Rotor leak. react. & $X_{lr} = 0.1104$ \\ 
        Base frequency & $\SI{50}{\hertz}$ & Main reactance & $X_m = 2.3489$ \\ 
		Rotational speed & $\SI{596}{rpm}$ & Number of pole pairs & $p = 5$ \\
		dc-link voltage & $\SI{1.5937}{\pu}$ & & \\
		\hline 
	\end{tabular}
\end{table}

\begin{table}[t]
	\centering
	\caption{Controller and simulation parameters.}
	\label{tab:control_parameters}
	\begin{tabular}{ccc}
		\hline
		Parameter name & Parameter symbol & Parameter values \\
		\hline
		torque weighting factor & $\lambda_\T$ & $0.052$ \\
		switching weight & $\lambda_\U$ & $3.8\times10^{-3}$\\
		controller sampling interval & $T_\C$ & $\SI{25}{\micro\second}$\\
		horizon & $N$ & $5$\\
		upper bound on parent nodes & $n_{\p,\max}$ & $\nodelim$\\
		simulation sampling interval & $T_\s$ & $\SI{.5}{\micro\second}$ \\
		\hline
	\end{tabular}
\end{table}

Physical parameters are provided in Table \ref{tab:sim_parameters} and they are taken from~\cite[Tab.~7.13]{geyer2005low}.
The control parameters of Table~\ref{tab:control_parameters} are chosen as follows. The torque weighting factor, $\lambda_\T = 0.052$, is adopted from \cite[$\S4.3$]{book_geyer}. The switching penalty $\lambda_\U = 13\times10^{-3}$ is tuned to achieve a device switching frequency of $f_\sw = \SI{215}{\hertz}$. The sampling interval is $T_\C = \SI{25}{\micro\second}$, and the prediction horizon is $N=5$.
Our comparisons will be based on the steady-state behaviour and on the step-up of the torque reference $T_\e^*$. The torque step-up will occur at $t=0.705$s from $T_\e^* = 0.2$ to $T_\e^* = 1$. During all studies, the stator flux magnitude reference will be kept at a constant level $\Psi_\s^* = 1$. The system is simulated with the time step of $T_\s = \SI{0.5}{\micro\second}$. A measurement noise of $\eta \sim \mathcal{U}(-2.5 \times 10^{-3},\ 2.5 \times 10^{-3})$ pu is imposed, where $\mathcal{U}(a,\ b)$ is the uniform distribution in the interval $[a,\ b]$.

The number of nodes traversed by Algorithm~\ref{alg:1} is assumed to be a good indicator of the complexity of the branch-and-bound. In fact, we computed the Pearson correlation coefficient of $\mathrm{Corr}(T_b,n_p) = 0.9989$ for the computational time $T_b$ and the number of nodes traversed $n_\p$. 
Considering this, as a node limit, $n_{\p,\max} = \nodelim$ is imposed in our simulation environment. On the other hand, the \ac{SDP} is solved via SCS~\cite{o2016conic} called through YALMIP~\cite{lofberg2004yalmip}, with an iteration limit of $120$. For the \ac{SDP} under consideration, this limit influences the high-order precision, which is later rounded during the extraction process.
As a remark, an assessment on whether these limits provide a fair comparison requires an implementation on an embedded platform, such as an FPGA.
The goal of this paper is to showcase the capabilities of \ac{SDP} relaxations for direct torque control and bring the two research communities closer together.
As a further justification of this choice, computation limits tighter than the one we imposed on the branch-and-bound have already been demonstrated for the sphere-decoder based on the number of operations, see~\cite{karamanakos2015suboptimal}.

\subsection{Efficiency analysis for branch-and-bound}

This study shows how the node limit $n_\p = \nodelim$ could be violated significantly by Algorithm~\ref{alg:1} if we do not impose it explicitly within the algorithm. Specifically, the efficiency of the branch-and-bound relies on two conditions: \textit{(i)} the accuracy of the initial guess, \textit{(ii)} the operational condition of the system. We study both in the following.

The impact of the first condition can be observed during the steady-state. We warm-start Algorithm~\ref{alg:1} with an all-zero initial guess, an educated guess, and the optimal solution.
We include the optimal solution here to specifically demonstrate how the effectiveness of the bounding decreases later during torque transients.
Table~\ref{tab:comp_complexity} shows the resulting maximum number of nodes traversed by the branch-and-bound. During steady-state, observe that the educated guess performs rather similar to the optimal solution, and both stay well under the limit $n_\p = \nodelim$. The all-zero initial guess, on the other hand, can result in violating the node limit $n_\p = \nodelim$, making it unsuitable for Algorithm~\ref{alg:1}. \looseness=-1
\begin{table}[t]
	\centering
	\caption{Number of nodes traversed for different initial guesses and situations. \textit{Max-steady-state} refers to the maximum during steady-state. \textit{Torque-step-instance} refers to the solution during the first instance of the step.}
	\label{tab:comp_complexity}
	\begin{tabular}{ccc}
		\hline
		Initial solution  & Max-steady-state & Torque-step-instance\\
		\hline
		Zero guess  & $1363$ & $69950$ \\
		Educated guess  & $361$ & $59674$ \\
		Optimal guess  & $358$ & $35524$ \\
		\hline
	\end{tabular}
\end{table}

The impact of operational condition on the efficiency is apparent during the torque step. Table~\ref{tab:comp_complexity} shows that, especially in the first instance after the step command, the number of traversed nodes peaks. Regardless of the initial guess, the node limit is significantly surpassed. This is because the inaccurate tracking during the step-up leads to much higher cost function values. The branch-and-bound algorithm, having no information on the future cost, needs to traverse deep into the tree before deciding on pruning the sub-trees. As an interesting remark, the educated guess performs almost as bad as the all-zero initial guess during transients.

\subsection{Performance comparison} \label{sec:numerica_performance}

\begin{figure}[!t]
    \centering
    \begin{subfigure}{.851\linewidth}
    \setlength{\fwidth}{\linewidth}
    \setlength{\fheight}{0.61\fwidth}
%
{\scriptsize
\begin{tikzpicture}
\def\dx{-0.7};
\def\dy{-0.4};
\begin{axis}[
grid = major,
xmin=0.65,
xmax=0.75,
xlabel = {\scriptsize Time [\si{\second}]},
ylabel = {\scriptsize $T_\e$ [\si{\pu}]},
ylabel near ticks,
xlabel near ticks,
width=\fwidth,
height=\fheight,
legend style={at={(0.05,1.08)},anchor=north west},
ymin = -0.7,
ymax = 1.2,
ticklabel style={font=\scriptsize},
label style={font=\scriptsize},
legend style={font=\scriptsize},
]
\addplot [restrict x to domain=0:2,color=black, thick, dotted] table[x=t, y=T_eref,mark=none] {Figures/torque.txt};
\addplot [restrict x to domain=0:2,color=red, thick, densely dashed] table[x=t, y=T_e2,mark=none] {Figures/torque.txt};
\addplot [restrict x to domain=0:2,color=blue, thick] table[x=t, y=T_e3,mark=none] {Figures/torque.txt};
\legend{reference, Algorithm~\ref{alg:1}, Algorithm~\ref{alg:2}}
\end{axis}
\end{tikzpicture}}
    \caption{Electromagnetic torque.}
    \label{fig:torque}
    \end{subfigure}
    \begin{subfigure}{.851\linewidth}
    \setlength{\fwidth}{\linewidth}
    \setlength{\fheight}{0.56\fwidth}
%
{\scriptsize
\begin{tikzpicture}
\def\dx{-0.7};
\def\dy{-0.4};
\begin{axis}[
grid = major,
xmin=0.65,
xmax=0.75,
xlabel = {\scriptsize Time [\si{\second}]},
ylabel = {\scriptsize $\Psi_\s$ [\si{\pu}]},
ylabel near ticks,
xlabel near ticks,
width=\fwidth,
height=\fheight,
legend style={at={(0.05,1.07)},anchor=north west},
ymin = .875,
ymax = 1.19,
ticklabel style={font=\scriptsize},
label style={font=\scriptsize},
legend style={font=\scriptsize},
]
\addplot [restrict x to domain=0:2,color=black, thick, dotted] table[x=t, y=Psi_sref,mark=none] {Figures/stator_flux.txt};
\addplot [restrict x to domain=0:2,color=red, thick, densely dashed] table[x=t, y=Psi_s2,mark=none] {Figures/stator_flux.txt};
\addplot [restrict x to domain=0:2,color=blue, thick] table[x=t, y=Psi_s3,mark=none] {Figures/stator_flux.txt};
\legend{reference, Algorithm~\ref{alg:1}, Algorithm~\ref{alg:2}}
\end{axis}
\end{tikzpicture}}
    \caption{Stator flux magnitude.}
    \label{fig:stator}
    \end{subfigure}
     \begin{subfigure}{.851\linewidth}
    \setlength{\fwidth}{\linewidth}
    \setlength{\fheight}{0.56\fwidth}
%
{\scriptsize
\begin{tikzpicture}
\def\dx{-0.7};
\def\dy{-0.4};
\begin{axis}[
grid = major,
xmin=0.65,
xmax=0.75,
xlabel = {\scriptsize Time [\si{\second}]},
ylabel = {\scriptsize $\b{i}_\s$ [\si{\pu}]},
ylabel near ticks,
xlabel near ticks,
width=\fwidth,
height=\fheight,
legend style={at={(0.05,1.13)},anchor=north west},
ymin = -1.2,
ymax = 1.2,
ticklabel style={font=\scriptsize},
label style={font=\scriptsize},
legend style={font=\scriptsize},
]

\addplot [restrict x to domain=0:2,color=red, thick, densely dashed] table[x=t, y=i_s2_1,mark=none] {Figures/currents.txt};
\addplot [restrict x to domain=0:2,color=blue, thick] table[x=t, y=i_s3_1,mark=none] {Figures/currents.txt};

\addplot [restrict x to domain=0:2,color=red, thick, densely dashed] table[x=t, y=i_s2_2,mark=none] {Figures/currents.txt};
\addplot [restrict x to domain=0:2,color=blue, thick] table[x=t, y=i_s3_2,mark=none] {Figures/currents.txt};

\addplot [restrict x to domain=0:2,color=red, thick, densely dashed] table[x=t, y=i_s2_3,mark=none] {Figures/currents.txt};
\addplot [restrict x to domain=0:2,color=blue, thick] table[x=t, y=i_s3_3,mark=none] {Figures/currents.txt};

\legend{Algorithm~\ref{alg:1}, Algorithm~\ref{alg:2}}
\end{axis}
\end{tikzpicture}}
    \caption{Three-phase stator current measurements.}
    \label{fig:currents}
    \end{subfigure}
    \caption{Torque reference step-up simulation results. Depicted are \ac{c1} with solvers as in Algorithms~\ref{alg:1} and~\ref{alg:2}.}
\end{figure}
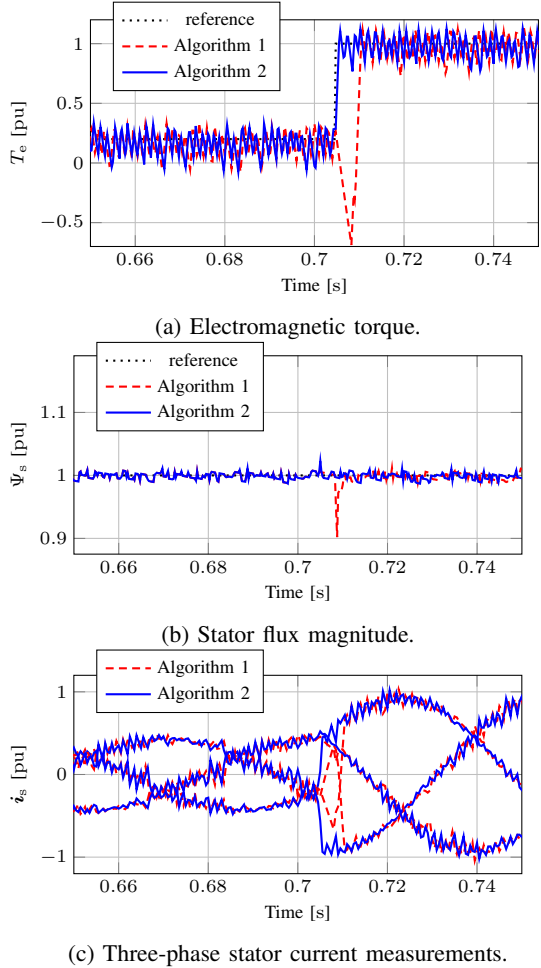

 Figure~\ref{fig:torque} shows the torque trajectories of \ac{c1} with Algorithms~\ref{alg:1} and~\ref{alg:2}. The input sequence from the node-limited Algorithm~\ref{alg:1} cannot track the torque reference, whereas Algorithm~\ref{alg:2} can find a close-to-optimal input sequence thanks to the SDP, leading to high dynamical performance. Figures~\ref{fig:stator} and~\ref{fig:currents} show that Algorithm~\ref{alg:1} would create large deviations in both stator flux and currents.
 As a remark, all the aforementioned input sequence extraction methods gave the same result during the steady-state in terms of harmonic distortions since branch-and-bound achieved optimality. Moreover, all these methods were capable of fixing the problematic instances observed during large torque transients due to branch-and-bound node limit.

 The \ac{SDP} relaxations have been demonstrated to be consistently and reliably solved within an acceptable time frame. As mentioned before, we used SCS~\cite{o2016conic} called through YALMIP~\cite{lofberg2004yalmip}. The average solver total time (includes both the setup and the solve stages of SCS) was \SI{8.3}{\milli\second} and the maximum total time was \SI{10.4}{\milli\second} for both the steady-state and the transient phases.\footnote{All problems are solved via MATLAB on a computer equipped with \SI{32}{\giga\byte} RAM and a \SI{2.5}{\giga\hertz} Intel i7 processor.} 
 Considering the times reported in~\cite[Fig. 5]{stellato_2017}, one could obtain more than two-order improvement in execution times by an FPGA. Nonetheless, this observation needs to be verified for the \ac{SDP}s on an embedded platform. \looseness=-1

\section{Conclusion}
We formulated the \ac{SDP} relaxation of \ac{c1}, and proposed solving it in parallel to a node-limited branch-and-bound designed for the original problem. Case studies showcased that a solution extracted from the \ac{SDP} can be advantageous during torque transients. Future work could solve the \ac{SDP} relaxation of \ac{c1} on an FPGA, and provide a theoretical approximation guarantee on the extracted input sequences.\looseness=-1

\bibliographystyle{IEEEtran}
\bibliography{references.bib}

\appendices
\section{Matrices for the \ac{SDP} formulation}\label{app:matrices}

The matrices $\b{\Gamma}_{\s,\ell}$, $\b{\Gamma}_{\r,\ell}$, $\b{\Upsilon}_{\s,\ell}$, $\b{\Upsilon}_{\r,\ell}$ are defined as follows
{\medmuskip=0.1mu\thinmuskip=0.1mu\begin{align*}
    \b{\Gamma}_{\s,\ell} &= [\b{A}_{11}\ \b{A}_{12}] \b{A}^{\ell-1},\quad\quad
    \b{\Upsilon}_{\s,1} = [\b{B}_1\ \b{0}_{2\times 3(N-1)}],\\
    \b{\Upsilon}_{\s,\ell} &= \Big[[\b{A}_{11}\ \b{A}_{12}]\b{A}^{\ell-2}\b{B}\ \ldots\ [\b{A}_{11}\ \b{A}_{12}]\b{B}\ \b{B}_1 \ \b{0}_{2\times 3(N-\ell)}\Big],\\
    &\hspace{6cm} \forall \ell>1,\\
    \b{\Gamma}_{\r,\ell} &= [\b{A}_{21}\ \b{A}_{22}] \b{A}^{\ell-1},\quad\quad
    \b{\Upsilon}_{\r,1} = [\b{0}_{2\times 3} \ \b{0}_{2\times 3(N-1)}],\\
    \b{\Upsilon}_{\r,\ell} &= \Big[[\b{A}_{21}\ \b{A}_{22}] \b{A}^{\ell-2}\b{B}\ \ldots\ [\b{A}_{21}\ \b{A}_{22}]\b{B}\ \b{0}_{2\times 3} \ \b{0}_{2\times 3(N-\ell)}\Big],\\
    &\hspace{6cm} \forall \ell>1,
\end{align*}}with $\b{A} = \begin{bmatrix}
    \b{A}_{11} & \b{A}_{12} \\
    \b{A}_{21} & \b{A}_{22}
\end{bmatrix}$ and $\b{B} = \begin{bmatrix}
    \b{B}_1 \\ \b{0}_{2\times3}.
\end{bmatrix}$ Here, $\b{I}_n \in \mathbb{R}^{n\times n}$ denotes the identity matrix and $\b{0}_{n\times m} \in \mathbb{R}^{n \times m}$ is an all-zero matrix of dimension $n$ times $m$.

Given these matrix definitions, torque can be obtained as 
\begin{align*}
    y_1(k+\ell) &= T_\factor \b{x}_{3:4}(k+\ell)^\top \b{J} \b{x}_{1:2}(k+\ell)\\
	&\hspace{-2cm}= T_\factor \left(\b{\Gamma}_{\r,\ell}\b{x}(k) + \b{\Upsilon}_{\r,\ell}\b{U}(k)\right)^\top \b{J} \left(\b{\Gamma}_{\s,\ell} \b{x}(k) + \b{\Upsilon}_{\s,\ell} \b{U}(k)\right)\\
    &\hspace{-2cm}=  \begin{bmatrix}
        1 & \b{U}(k)^\top
    \end{bmatrix} \b{Q}_\ell(k) \begin{bmatrix}
	    1 \\ \b{U}(k)
	\end{bmatrix},
\end{align*}
where
\begin{equation*}
    \b{Q}_\ell(k) =  T_\factor\begin{bmatrix}
		\b{x}(k)^\top \b{\Gamma}_{\r,\ell}^\top \b{J}\b{\Gamma}_{\s,\ell} \b{x}(k) & \b{x}(k)^\top\b{\Gamma}_{\r,\ell}^\top \b{J} \b{\Upsilon}_{\s,\ell} \\
		\b{\Upsilon}_{\r,\ell}^\top \b{J} \b{\Gamma}_{\s,\ell} \b{x}(k) & \b{\Upsilon}_{\r,\ell}^\top \b{J} \b{\Upsilon}_{\s,\ell}
	\end{bmatrix}.
\end{equation*}
Similarly, stator flux magnitude is given by
\begin{align*}
	y_2(k+\ell) &= ||\b{x}_{1:2}(k+\ell)||_2^2 \\
	&= \left(\b{\Gamma}_{\s,\ell} \b{x}(k) + \b{\Upsilon}_{\s,\ell} \b{U}(k)\right)^\top\left(\b{\Gamma}_{\s,\ell} \b{x}(k) + \b{\Upsilon}_{\s,\ell} \b{U}(k)\right)\\
	&= \begin{bmatrix}
		1 & \b{U}(k)^\top
	\end{bmatrix} \b{W}_\ell(k) \begin{bmatrix}
		1 \\
		\b{U}(k)
	\end{bmatrix},
\end{align*}
where
\begin{equation*}
    \b{W}_\ell(k) = \begin{bmatrix}
		\b{x}(k)^\top \b{\Gamma}_{\s,\ell}^\top\b{\Gamma}_{\s,\ell} \b{x}(k) & \b{x}(k)^\top\b{\Gamma}_{\s,\ell}^\top\b{\Upsilon}_{\s,\ell} \\
		\b{\Upsilon}_{\s,\ell}^\top\b{\Gamma}_{\s,\ell}\b{x}(k) & \b{\Upsilon}_{\s,\ell}^\top\b{\Upsilon}_{\s,\ell}
	\end{bmatrix}.
\end{equation*}
Finally, for switching transitions, define $\b{Z}_\A = [1\ 0\ 0]$, $\b{Z}_\B = [0\ 1\ 0]$, and $\b{Z}_\C= [0\ 0\ 1]$ for each phase.
For phase a, and for $\ell = 1$, we have
\begin{equation*}
	\Delta u_\A(k) = u_\A(k) - u_\A(k-1) = \b{Z}_{\A,1}(k)
	\begin{bmatrix}
		1 \\
		\b{U}(k)
	\end{bmatrix}, \nonumber
\end{equation*}
where
\begin{equation*}
    \b{Z}_{\A,1}(k) = \begin{bmatrix}
		-u_\A(k-1) & \b{Z}_\A & \b{0}_{1\times 3(N-1)}
	\end{bmatrix}.
\end{equation*}
On the other hand, for $\ell > 1$:
\begin{align*}
	\Delta u_\A(k+\ell-1) &= u_\A(k+\ell-1) - u_\A(k+\ell-2) \\
	&= \b{Z}_{\A,\ell}(k)
	\begin{bmatrix}
		1 \\
		\b{U}(k)
	\end{bmatrix}, \nonumber
\end{align*}
where
\begin{equation*}
    \b{Z}_{\A,\ell}(k) = \begin{bmatrix}
		0 & \b{0}_{1\times3(\ell-2)} & -\b{Z}_a & \b{Z}_a & \b{0}_{1\times3(N-\ell)}
	\end{bmatrix}.
\end{equation*}
Other phases can be defined similarly.

\end{document}